\documentclass[12pt,reqno,conference,bibother]{amsart} 
\usepackage{amssymb,latexsym,verbatim}
\usepackage[colorlinks]{hyperref}

\allowdisplaybreaks

\newcommand\Tr{\text{Tr}}
\newcommand{\A}{\mathcal A}
\newcommand\pair[2]{\langle#1,#2\rangle}
\renewcommand{\phi}{\varphi}
\newtheorem*{main}{Main Theorem}
\newtheorem{theorem}{Theorem}

\newtheorem{proposition}[theorem]{Proposition}
\newtheorem{cor}[theorem]{Corollary}

\theoremstyle{definition}

\newtheorem{rem}[theorem]{Remark}

\begin{document}

\title[Index sets in model theory]%
{The complexity of the index sets of $\aleph_0$-categorical theories and of
Ehrenfeucht theories}
\author[Lempp]{Steffen Lempp}
\address{Department of Mathematics\\
  University of Wisconsin\\
  Madison, WI 53706-1388\\
  USA}
\email{lempp@math.wisc.edu}
\author[Slaman]{Theodore A. Slaman}
\address{Department of Mathematics\\
  University of California\\
  Berkeley, CA 94720-3840\\
  USA}
\email{slaman@math.berkeley.edu}
\thanks{The first author's research was partially supported by NSF grant
DMS-0140120 and by the Deutsche Forschungsgemeinschaft through a Mercator
Guest Professorship. The second author's research was partially supported by
the Alexander von Humboldt Foundation and by NSF Grant~DMS-9988644. Both
authors also wish to thank the University of Heidelberg and Klaus Ambos-Spies
for their hospitality during the time this research was carried out.}

\begin{abstract}
We classify the computability-theoretic complexity of two index sets of
classes of first-order theories: We show that the property of being an
$\aleph_0$-categorical theory is $\Pi^0_3$-complete; and the property of
being an Ehrenfeucht theory $\Pi^1_1$-complete. We also show that the
property of having continuum many models is $\Sigma^1_1$-hard. Finally, as a
corollary, we note that the properties of having only decidable models, and
of having only computable models, are both $\Pi^1_1$-complete.
\end{abstract}
\maketitle

\section{The Main Theorem}\label{sec:main}
Measuring the complexity of mathematical notions is one of the main tasks of
mathematical logic. Two of the main tools to classify complexity are provided
by Kleene's arithmetical and analytical hierarchy. These two hierarchies
provide convenient ways to determine the exact complexity of properties by
various notions of \emph{completeness}, and to give lower bounds on the
complexity by various notions of \emph{hardness}. (See, e.g.,
Kleene~\cite{Kl52}, Soare~\cite{So87} or Odifreddi~\cite{Od89,Od99} for the
definitions.)

This paper will investigate the complexity of properties of a first-order
theory, more precisely, the complexity of a countable first-order theory
having a certain number of models. Recall that a theory is called
\emph{$\aleph_0$-categorical} if it has only one countable model up to
isomorphism, and an \emph{Ehrenfeucht theory} if it has more than one but
only finitely many countable models up to isomorphism. In order to measure
the complexity of these notions, we will use \emph{decidable} first-order
theories, i.e., sets $T$ of first-order sentences closed under inference such
that membership in $T$ can be determined effectively.

The principal result of this paper is now the following

\begin{main}
We classify the complexity of some properties of decidable first-order
theories as follows:
\begin{enumerate}
\item The property ``$T$ is an $\aleph_0$-categorical theory'' is
  $\Pi^0_3$-complete.
\item The property ``$T$ is an Ehrenfeucht theory'' is $\Pi^1_1$-complete.
\item The property ``$T$ has continuum many pairwise nonisomorphic models''
is $\Sigma^1_1$-hard.
\end{enumerate}
\end{main}

In this paper, a \emph{theory} will be a set of first-order sentences closed
under inference (so that $T \vdash \sigma$ iff $\sigma \in T$ for any
sentence $\sigma$).

The rest of this paper is devoted to the proof of our Main Theorem. In
section~\ref{sec:cat}, we will prove clause (1) of our Main Theorem. In
section~\ref{sec:Pi11}, we will show that the property of a theory being an
Ehrenfeucht theory is $\Pi^1_1$, giving the upper bound for clause (2) of our
Main Theorem. Finally, in section~\ref{sec:reduc}, we will prove a
simultaneous reduction
$$
(\Pi^1_1,\Sigma^1_1) \le_m (\text{Ehrenfeucht},\text{Continuum})
$$
(where Ehrenfeucht and Continuum are (the index sets of) the properties of
being an Ehrenfeucht theory or a theory having continuum many models,
respectively). This reduction gives the lower bound for clause (2) and proves
clause (3) of our Main Theorem.

\section{\texorpdfstring{The Proof for $\aleph_0$-Categoricity}%
{The Proof for Countable Categoricity}}\label{sec:cat} We first show that
``$T$'s being $\aleph_0$-categorical'' is $\Pi^0_3$ to obtain the upper
bound: A theory $T$ (with a characteristic function given by a partial
computable function $\phi_e$) is $\aleph_0$-categorical iff, by the
Ryll-Nardzewski Theorem,
\begin{gather*}
\text{$\phi_e$ is a total function, $T$ is a complete consistent theory and}\\
\forall n\, \text{($T$ has only finitely many $n$-types),}
\end{gather*}
i.~e., iff
\begin{gather*}
e \in \text{Tot}
 \text{ and } \forall \sigma\, (T \vdash \sigma \text{ implies }
    \sigma \in T) \text{ and }\\
 \forall \sigma\, (T \not\vdash \sigma\wedge\neg\sigma) \text{ and }
 \forall \sigma\, (T \vdash \sigma \text{ or } T \vdash \neg\sigma)
    \text{ and }\\
\forall n\, \exists m\, \exists \sigma_0,\dots,\sigma_m\,
 [\forall i \le m\, (T \not\vdash \sigma_i) \text{ and } \qquad\\ \qquad
  \forall i \le m\, \forall \tau\, (T \vdash \sigma_i \rightarrow \tau
    \text{ or }
                                    T \vdash \sigma_i \rightarrow \neg\tau)].
\end{gather*}
(where $\sigma$ ranges over all first-order sentences, and $\sigma_i$ and
$\tau$ range over all first-order formulas in $n$ variables $\overline x$
which we suppress above, respectively). Now, since ``$T \vdash \sigma$'' is
$\Delta^0_1$ for decidable theories, inspection shows that ``$T$'s being
$\aleph_0$-categorical'' is $\Pi^0_3$. (Recall here that we assume that $T$
is closed under inference.)

In order to show that ``$T$'s being $\aleph_0$-categorical'' is
$\Pi^0_3$-complete, we present the following construction of a decidable
theory $T_e$ (uniformly in an index $e$) such that
$$
\text{$\forall n\,(W_{f(n)}$ is finite) iff $T_e$ is $\aleph_0$-categorical},
$$
where $f=f_e$ is a computable function such that
$$
\forall n\, (W_{f(n)} \text{ is finite})
$$
is a $\Pi^0_3$-complete predicate (see, e.g., Soare \cite[p.~68]{So87}).
Without of loss of generality, we may assume that for every stage $s$, there
is exactly one pair $\pair n x$ such that $x$ enters $W_{f(n)}$ at stage~$s$.
We denote this $n$ by $n_s$.

The signature of our theory $T_e$ now consists of relation symbols $R^n_s$
(for all $n,s \in \omega$) where each $R^n_s$ is an $n$-ary relation symbol.
At stage~$s$, we completely specify the relations $R^n_s$ for all $n \in
\omega$ as follows: For all $n \not= n_s$, we let the relation $R^n_s$ be
empty. For $n = n_s$, we let the relation $R^n_s$ be ``random'' over all
$R^{n_{s'}}_{s'}$ (for $s'<s$) in the sense that all finite extensions
consistent with the theory enumerated before stage~$s$ (in the relation
symbols $R^{n_{s'}}_{s'}$ for all $n$ and all $s' < s$) are realized by the
relation symbol $R^{n_s}_s$ added at stage~$s$.

To verify that the construction yields the theory $T_e$ with the desired
properties, first assume that for some $n$, $W_{f(n)}$ is infinite. Then
$R^n_s$ is nonempty for infinitely many $s$, and in fact the reduct of $T_e$
to these $n$-ary relations has continuum many consistent $n$-types, making
$T_e$ not $\aleph_0$-categorical.

On the other hand, assume that $W_{f(n)}$ is finite for all $n$. Then the
$n$-type of any $n$-tuple $\overline x$ is determined by the finitely
many nonempty relations of arity $\le n$ satisfied by $\overline x$, so
there are only finitely many $n$-types, and $\aleph_0$-categoricity
follows.

In fact, it is not hard to see that the theory $T_e$ admits elimination
of quantifiers (effectively uniformly in $e$) and is decidable (again
uniformly in $e$).

\section{The Upper Bound for Ehrenfeucht theories}\label{sec:Pi11}

\begin{proposition}\label{Pi11}
The set $\{e:\text{$T_e$ is an Ehrenfeucht theory}\}$ is $\Pi^1_1$.
\end{proposition}

\begin{proof}
Proposition~\ref{Pi11} follows from Sacks's \cite{Sa83} proof that every
countable model of an Ehrenfeucht theory has a hyperarithmetic
representation.

Suppose that $T$ is a recursive, complete, first-order, Ehrenfeucht theory.
From Sacks's theorem, each of the finitely many isomorphism types of models
of $T$ has a representative which is hyperarithmetically coded. In the same
paper, Sacks shows that any model $\A$ of $T$ completes its canonical Scott
analysis in finitely many steps. Consequently, the canonical Scott sentence
for $\A$ is uniformly hyperarithmetic in any representation $A$ of $\A$.
Further, if $A$ and $H$ are representations of countable first-order
structures and have the same finite-rank canonical Scott sentence, then there
is an isomorphism $\pi$ between the models coded by $A$ and $H$ which is
uniformly hyperarithmetically definable from $A$ and $H$. See, for example,
\cite{Na74} in which Nadel analyzes of Scott sentences and isomorphisms
between countable models within admissible sets. Thus, for $T$ with the given
properties, there are finitely many hyperarithmetic presentations
$H_1,\dots,H_k$ of models of $T$ such that for every representation $A$ of a
model of $T$, there is an isomorphism between the models coded by $A$ one of
the $H_i$'s which is hyperarithmetic relative to $A$ and $H_1,\dots,H_k$.

Conversely, for any recursive, complete, first-order theory $T$, if there is
a hyperarithmetic finite sequence $H_1,\dots,H_k$ of representations of
models of $T$ such that for every $A$, if $A$ codes a model of $T$ then there
is an isomorphism $\pi$ between the model coded by $A$ and that coded by one
of the $H_i$'s, then $T$ is an Ehrenfeucht theory.

Now, the proposition follows. $T_e$'s being a complete first-order theory is
an arithmetic property and hence $\Pi^1_1$. By the Spector-Gandy Theorem, see
Sacks~\cite{Sa90}, the $\Pi^1_1$ predicates are closed under existential
quantification over the hyperarithmetic sets. Thus, the condition above, that
there exist hyperarithmetic sets $H_1,\dots,H_k$, for all sets $A$, there
exists $\pi$ hyperarithmetic in $H_1,\dots,H_k$ and $A$, with arithmetically
described properties, is a $\Pi^1_1$ condition.
\end{proof}

\section{A Simultaneous Reduction for Ehrenfeucht Theories and Theories with
Continuum Many Models}\label{sec:reduc} In this section, we will establish
the simultaneous reduction
\begin{equation}\label{eq:reduc}
(\Pi^1_1,\Sigma^1_1) \le_m (\text{3Models},\text{Continuum})
\end{equation}
where 3Models and Continuum are (the index sets of) the properties of being
an Ehrenfeucht theory with exactly three countable models, and a theory
having continuum many countable models, respectively. Our proof will be based
largely on Reed~\cite{Re91}, which in turn used previous work of
Peretyat'kin~\cite{Pe73} and Millar~\cite{Mi83}. (Since our proof uses the
rather involved machinery of Reed~\cite{Re91} so heavily and mostly without
changes, we will assume familiarity with this paper throughout the rest of
the proof.)

We first observe that one can easily modify the proof of the
$\Sigma^1_1$-completeness of the property of a computable tree $\Tr \subseteq
\omega^{<\omega}$ having an infinite path to obtain a reduction
\begin{equation}\label{eq:reduc1}
(\Pi^1_1,\Sigma^1_1) \le_m (\text{NoPath},\text{InfPath})
\end{equation}
where NoPath and InfPath are (the index sets of) the properties of being a
computable tree $\Tr \subseteq \omega^{<\omega}$ having no infinite path, and
having continuum many infinite paths, respectively.

We now compose the reduction from (\ref{eq:reduc1}) with the reduction given
in Reed~\cite{Re91}: Reed defines, for each computable tree $\Tr \subseteq
\omega^{<\omega}$, and uniformly in an index $e$ of $\Tr$, a complete
decidable theory $T_e$ ``coding'' the tree $\Tr$ into a ``dense tree''.
(Actually, Reed only defines $T_e$ for trees $\Tr$ having exactly one
infinite path, but his definition in Part II of \cite{Re91} can be applied to
any computable $\Tr \subseteq \omega^{<\omega}$ and always yields a complete
decidable theory $T_e$.)

Checking over Reed's analysis of 1-types over $T_e$ in Part III of
\cite{Re91}, one can now easily verify the following: If $\Tr$ has no
infinite path then $T_e$ has exactly three countable models, namely, the
countable computable models omitting the type $\Gamma^*(x)$ (since $\Tr$
contains no infinite paths). On the other hand, if $\Tr$ has infinitely many
infinite paths then $T_e$ admits a partial type $\Gamma^*_f = \{c_\xi<x \mid
\xi \subset f\}$ for each infinite path $f \in [\Tr]$, and by an argument
analogous to that for Corollary 15.1 in Reed~\cite{Re91} (which shows that
the type $\Delta(x)$ can be realized without the type $\Gamma^*(x)$ being
realized), for any two distinct paths $f, g \in [\Tr]$, the partial types
$\Gamma^*_f$ and $\Gamma^*_g$ can be realized independently of each other.
Thus $T_e$ has continuum many types and so also continuum many countable
models.

To finish off, we note an easy corollary of our proof as well as two remarks:

\begin{cor}
The property of a decidable theory having only decidable models is
$\Pi^1_1$-hard, as is the property of a decidable theory having only
computable models.
\end{cor}

\begin{rem}
\begin{enumerate}
\item Our proof actually shows that the property of being an Ehrenfeucht
theory with exactly three models is already $\Pi^1_1$-complete.
\item The first-order language used in our proof is not fixed; i.~e., the
language depends on the tree $\Tr$ and thus on the index $e$ used for the
construction. However, we can simply add constant symbols $c_\eta$ (for $\eta
\in \omega^{<\omega} - T$) denoting a fixed element of our model, and empty
relations $E^\eta_\xi$, etc.\ (for $\eta$ or $\xi \in \omega^{<\omega} - T$),
to achieve a fixed language.
\end{enumerate}
\end{rem}

\end{document}